\documentclass[12pt]{amsart}
\usepackage{amssymb,amsmath,graphicx,verbatim,amsthm}
\usepackage{color}
\usepackage{xcolor}
\usepackage{mdframed}

\definecolor{env_back}{gray}{0.8}
\definecolor{thm_color}{rgb}{0,0,1}

\newtheorem{thm}{Theorem}
\newtheorem{cor}[thm]{Corollary}
\newtheorem{lem}[thm]{Lemma}

\newtheorem{prop}[thm]{Proposition}

\newtheorem{conj}[thm]{Conjecture}

\newtheorem*{clm*}{Claim}

\theoremstyle{definition}
\newtheorem{dfn}[thm]{Definition}
\newtheorem{exm1}[thm]{Example}

\theoremstyle{remark}
\newtheorem{rem}[thm]{Remark}

\newenvironment{lem*}[1]{\vspace{1ex}\noindent
{\bf Lemma* (#1).} [restatement]  \hspace{0.5em} \em }{ }
\newenvironment{thm*}[1]{\vspace{1ex}\noindent 
{\bf Theorem* (#1).} [restatement]  \hspace{0.5em} \em }{ }

{\begin{exm1}}%
{\hfill$\bigtriangleup\bigtriangledown\bigtriangleup$
\end{exm1} }

\newcommand{\set}[1]{\left\{#1\right\}}

\newcommand{\N}{\mathbb{N}}

\newcommand{\R}{\mathbb{R}}

\newcommand{\eps}{\varepsilon}

\DeclareMathOperator{\E}{\mathbb{E}}    
\renewcommand{\Pr}{}
\let\Pr\relax
\DeclareMathOperator{\Pr}{\mathbb{P}}

\newcommand{\1}[1]{\mathbf{1}_{\set{ #1 } }}

\author{Idan Perl}

\usepackage{float} 

\newcommand{\F}{\hat{F}}

\title{SIT measures and Transience}

\begin{document}

\begin{abstract}
A graph is called (edge)-SIT if for some probability measure on paths, the number of mutual edges in two independent paths has finite mean. We show that transience is equivalent to the property of (edge)-SIT.
\end{abstract}  

\maketitle

\section*{Introduction}
Let $G$ be a connected graph and let $\mu$ be a probability measure on paths in $G$. 
In this note we consider the intersections of two independent paths of law $\mu$.
The motivation for this originated in the paper \cite{BPP97}. There, the authors show that for any graph $G$, if there exist a measure $\mu$ such that this intersection has an exponentially decaying tail, 
then oriented percolation clusters are transient - as long as the retention parameter $p$ is close enough to $1$. 
In this note we consider graphs which admit a measure $\mu$ such that the intersection size  
is an \textit{integrable} random variable. 
We inspect the relationship between this property and transience of the graph. 

We start by formally defining the above concepts.
Let $G$ be a connected graph and fix a root vertex $o$. Denote by $\Gamma_o$ the set of simple paths in $G$ starting at $o$:
$$
\Gamma_o:=\{ \gamma:\N \hookrightarrow G:\ \gamma_0=o,\ \gamma_{n+1}\sim \gamma_n  \}.
$$
We say that a directed edge $e=(x,y)$ is in a path $\gamma$, denoted $(x,y) \in \gamma$, 
if there exists $n\in\N$ such that $\gamma_n=x,\gamma_{n+1}=y$. 
Similarly, a vertex $x$ is in $\gamma$, denoted $x \in \gamma$, if there exists $n$ such that $\gamma_n= x$.
For a probability measure $\mu$ on $\Gamma_o$, denote by $\mu\otimes\mu$ the canonical measure on $\Gamma_o\times\Gamma_o$ (of two independent paths of law $\mu$). 
For two paths $\alpha,\beta\in\Gamma_o$, let 
\begin{align}
\label{dfn: vertex intersection}
|\alpha\cap\beta|_V := \sum_{v\in G} \1{v\in\alpha \ , \ v \in \beta}
\end{align}
be the vertex intersection of the paths.
Similarly, let $|\alpha\cap\beta|_E$ be the edge intersection of $\alpha$ and $\beta$:
\begin{align}
\label{dfn: vertex intersection}
|\alpha\cap\beta|_E := \sum_{\substack{u,v\in G \\ u\sim v}} \1{(u,v)\in\alpha \ , \ (u,v) \in \beta}.
\end{align}
\begin{dfn}
A probability measure $\mu$ on $\Gamma_o$ is called \textit{vertex}-SIT (summable intersection tail) if
\begin{align*}
\E_{\mu\otimes\mu}|\alpha\cap\beta|_V<\infty,
\end{align*}
and \textit{edge}-SIT if
\begin{align*}
\E_{\mu\otimes\mu}|\alpha\cap\beta|_E<\infty.
\end{align*}
\end{dfn}
Note that if $\mu$ is vertex-SIT then it is edge-SIT. The converse is not true in general.

\begin{dfn}
A graph $G$ is called vertex-SIT (resp. edge-SIT) if there exists a vertex $o$ and a vertex-SIT (resp. edge-SIT) probability measure on $\Gamma_o$. 
\end{dfn}




We will prove the following theorem:
\begin{thm}
\label{thm: sit iff transient}
Let $G$ be a graph. $G$ is transient if and only if it is edge-SIT. If furthermore $G$ has uniformly bounded degree, another equivalent condition is that $G$ is vertex-SIT. 
\end{thm}

Following the proof of Theorem \ref{thm: sit iff transient}, we give an example of a transient unbounded degree graph which is not vertex-SIT, implying that bounded degree condition - while not necessary - is not redundant.

Before moving to the proof let us remark regarding the results of \cite{BPP97}.
There the notion of EIT (exponential intersection tails) was introduced.
A probability measure $\mu$ on $\Gamma_o$ is called EIT (summable intersection tail) if
there exists $\eps>0$ such that 
\begin{align*}
\E_{\mu\otimes\mu} e^{\eps |\alpha\cap\beta|_E } <\infty .
\end{align*}
$G$ is said to be EIT if there exists an EIT measure on $\Gamma_o$.
In \cite{BPP97} it is shown that if $G$ is EIT then for some large enough parameter $0< p < 1$,
Bernoulli percolation on $G$ admits transient infinite components.

The following has been conjectured by A.\ Yadin.
\begin{conj}
Let $G$ be a Cayley graph of a finitely generated group.
Then, $G$ is SIT if and only if $G$ is EIT.
\end{conj}

For more on the connection between transience of a group and EIT see \cite{RY} and references therein.

\section*{Proof}

We begin with the simpler implication.
\begin{prop}
If $G$ is edge-SIT then it is transient.
\end{prop}
\begin{proof}
Let $\mu$ be an edge-SIT measure on $\Gamma_o$. We will use it to construct a finite energy positive flow on $G$.  For a directed edge $e=(v,u)$, define 
\begin{align*}
F(v,u):=\Pr_{\mu}[(v,u)\in\gamma]-\Pr_{\mu}[(u,v)\in\gamma].   
\end{align*}
It is immediate that $F$ is anti-symmetric. Let $v\neq o$. Since paths are infinite, a path $\gamma$ visits $v$ if and only if it leaves it. That is, $\sum_{w\sim v} \1{(v,w)\in\gamma}-\1{(w,v)\in\gamma}=0$. Since $F(v,u)=\E_{\mu}[\1{(v,u)\in\gamma}-\1{(u,v)\in\gamma}]$, we deduce
\begin{align*}
\mathrm{div} F(v)=\sum_{u\sim v} F(v,u)=\E_{\mu}\big [\sum_{w\sim v} \1{(v,w)\in\gamma}-\1{(w,v)\in\gamma} \big]
=0.
\end{align*}
For the root $o$, since paths are simple, we have $\Pr_{\mu}[(u,o)\in\gamma]=0$ for all $u\sim o$. Also, $\Pr_{\mu}[(o,u)\in\gamma]>0$ for some $u\sim o$, hence $\mathrm{div} F(v)\neq 0$. Finally, we verify that $F$ has finite energy:
\begin{align*}
\sum_{v\sim u} F(v,u)^2
&=\sum_{v\sim u}\big(\Pr[(v,u)\in\gamma]-\Pr[(u,v)\in\gamma]\big)^2	\\
&\leq 2\cdot \sum_{v\sim u} \Pr[(v,u)\in\gamma]^2	\\
&=2\cdot \E_{\mu\otimes\mu}|\alpha\cap\beta|_E
<\infty
\end{align*}
since $\mu$ is edge-SIT.
\end{proof}


Before proving the implication of Theorem \ref{thm: sit iff transient}, we write some preliminaries. Let $G$ be a transient graph, and fix a root $o$. Let $F$ be a unit flow on $G$ from $o$ to $\infty$. For a vertex $v$, denote by $f(v)$ the total positive flow exiting $v$:
\begin{align}
\label{eqn: f}
f(v):
=\sum_{\substack{ w\sim v \\ \F(v,w)>0}}F(v,w).
\end{align} 
For $v\neq 0$, since $\mathrm{div} F(v)=0$, this is also equal to the amount of flow entering $v$:
\begin{align}
\label{eqn: f again}
v\neq 0 \implies f(v)
=\sum_{\substack{ w\sim v \\ F(w,v)>0}}F(w,v).
\end{align}
For $u,v\in G$, let   	
$$
Q(u,v):=f(u)^{-1}\cdot F(u,v)\cdot \1{F(u,v)>0}
$$
if $f(u)>0$, and $Q(u,v):=0$ otherwise. This defines a transition kernel on $G$, and we denote by $\mu_F$ the corresponding probability measure on random walks starting at $o$. Seen as an operator, $Q$ acts on functions $\varphi:G\to\R$ from the right by 
$$
(\varphi Q)(v):=\sum_{u\in G} \varphi(u)\cdot Q(u,v).
$$

The following lemma implies that we can choose $F$ so that $\mu_F$ will be supported on simple paths. We say that $\gamma=(v_0,v_1,...,v_n=v_0)$ is a \textit{positive energy loop} if $F(v_i,v_{i+1})>0$ for all $0\leq i\leq n-1$.

\begin{lem}
\label{lem: no loops}
Suppose $G$ is a transient graph. There exists a unit flow $F$ that admits no positive flow loops. Moreover, $F$ satisfies $F(u,o)\leq 0$ for all $u\in G$. 
\end{lem}
\begin{proof}
For a vertex $x$, let $\Pr_x$ be the probability measure on simple random walk paths starting at $x$. Denote by $T_o$ the hitting time of the distinct vertex $o$. By a classic construction, the following function is a finite energy flow on $G$:
\begin{align*}
F(x,y):=\Pr_x[T_o<\infty]-\Pr_y[T_o<\infty].
\end{align*}
It is straightforward to check that this function induces no positive energy loops and that $F(u,o)\leq 0$ for all $u\in G$.
\end{proof}

\begin{lem}
\label{f-g}
Suppose $G$ is a transient graph, and let $F$ be a flow as in Lemma \ref{lem: no loops}. Let $\mu_F$ be the corresponding measure on $\Gamma_o$. Let $g(v)=\Pr_{\mu_F}[T_v<\infty]$, the probability that a $\mu_F$-generated path will visit $v$. Let $f$  be as in \eqref{eqn: f}. Then for every $v\in G$, $g(v)\leq f(v)$.
\end{lem}  
\begin{proof}
Let $v\neq o$. By \eqref{eqn: f again}, we have
\begin{align*}
(fQ)(v)&=\sum_{u\in G} f(u)\cdot Q(u,v)
	    =\sum_{\substack{u\in G \\ f(u)>0}} 
	    f(u)\cdot f(u)^{-1}\cdot F(u,v)\cdot \1{F(u,v)>0}\\
	   &=f(v)
\end{align*}

By Lemma \ref{lem: no loops}, $Q(u,o)=0$ for every $u\in G$, hence $(fQ)(o)=0$. Also, since $F$ is a unit flow we have $f(o)=1$. Together, denoting by $Q^0$ the identity matrix, we get $fQ=f-\delta_0=f-Q^0(o,\cdot)$. Using induction, we get that for all $v$ and $n$, 
\begin{align}
(fQ^n)(v)=f(v)-\sum_{k=0}^{n-1}Q^k(o,v).
\end{align}
Moreover, the entries of $f$ and $Q$ are positive, thus so are the entries of $fQ^n$ for all $n$. Now, $Q^k(u,v)$ is the probability that a $\mu_F$-generated path $\gamma$ will get from $u$ to $v$ in $k$ steps. Since - by Lemma \ref{lem: no loops} - $\mu_F$ is supported on simple paths, we deduce $Q^k(o,v)=\Pr_{\mu_F}[T_v=k]$. Hence,  
\begin{align}
\label{g}
g(v)=\Pr_{\mu_F}[T_v<\infty]=\sum_{k=0}^{\infty} \Pr_{\mu_F}[T_v=k]=
\sum_{k=0}^{\infty} Q^k(o,v). 
\end{align}


Altogether, for all $v$, we get
\begin{align*}
0\leq \lim_{n\to\infty} (fQ^n)(v)
 =	  \lim_{n\to\infty} f(v)-\sum_{k=0}^{n-1}Q^k(0,v) 	
 = f(v)-g(v)
\end{align*}
as required.
\end{proof}

We are now ready to show:

\begin{prop}
\label{hard direction}
Let $G$ be a graph. If $G$ is transient then $G$ is edge-SIT. If furthermore $G$ has uniformly bounded degree, then $G$ is vertex-SIT. 
\end{prop}

\begin{proof}

Let $G$ be a transient graph. Let $F$ be a unit flow from $o$, and $\mu_F$ the corresponding measure. Let $e=(u,v)$ be a directed edge with $F(e)>0$. By $\eqref{g}$ and Lemma \ref{f-g}, we have
\begin{align*}
\Pr_{\mu_F}[(u,v)\in\gamma]
&=\sum_{k=0}^{\infty}Q^k(o,u)Q(u,v)
=g(u)\cdot Q(u,v)	\\
&\leq f(u)\cdot Q(u,v)=F(e),
\end{align*}
and we deduce
\begin{align*}
\E_{\mu_F \otimes \mu_F} |\alpha\cap\beta|_E=\sum_{F(e)>0}\Pr[e\in\gamma]^2\leq \sum_{e\in E} F(e)^2<\infty.
\end{align*}

For the second part of the proposition, by Cauchy-Schwarz inequality we have
\begin{align*}
f(v)^2
&=\big( \sum_{\substack{u\sim v \\ F(u,v)>0}}F(u,v) \ \big)^2	\\
&\leq deg(v) \cdot \sum_{\substack{u\sim v \\  F(u,v)>0}}F(u,v)^2	
\end{align*}
Hence, by Lemma \ref{f-g}, if there is a constant $d$ such that $deg(v)\leq d$, we have
\begin{align*}
\E_{\mu_F \otimes \mu_F} |\gamma_1\cap\gamma_2|_V
&=\sum_{v\in G}\Pr[v\in \gamma]^2=\sum_{v\in G} g(v)^2 \\
&\leq \sum_{v\in G} f(v)^2\leq d\cdot \sum_{u\sim v}F(u,v)^2<\infty,
\end{align*}
as required. 
\end{proof}

\begin{rem}
As can be seen in the proof of proposition \ref{hard direction}, bounded degree is not necessary. It is enough to have a degree function that satisfies
\begin{align*}
\sum_{u\sim v}deg(v)\cdot F(u,v)^2<\infty.
\end{align*}
\end{rem}

We proceed to show an example of a transient graph which is edge-SIT but not vertex-SIT. Consider the following finite graph:
\begin{figure}[H]
\includegraphics[scale=.7]{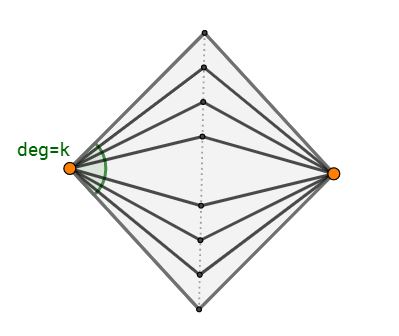}
\end{figure}

That is, two vertices mediated by $k$ (degree 2) vertices. We call this graph a degree $k$ diamond. Now, consider the standard graph of $\N$: vertices $m,n$ are connected with an edge if $|m-n|=1$. Replace every edge $(n-1,n)$ with a degree $2^n$ diamond. Call this graph $G$.

\begin{figure}[H]
\includegraphics[scale=.7]{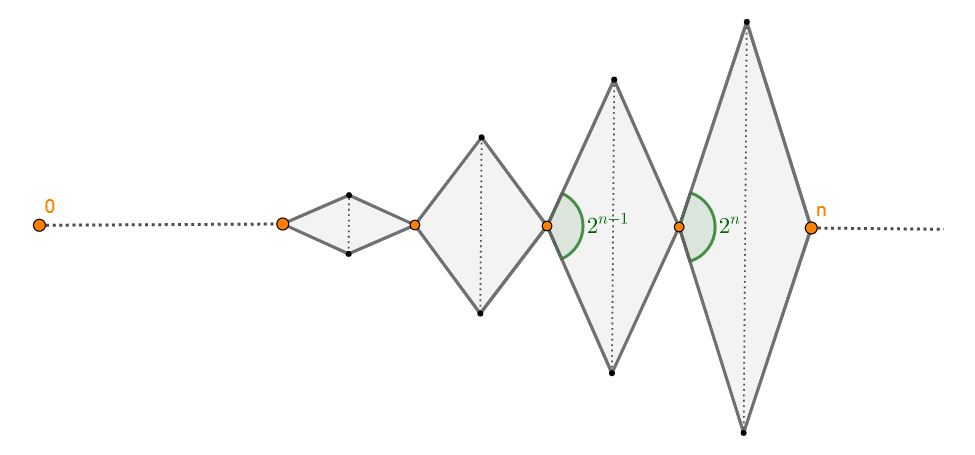}
\end{figure}

By network reduction, the macroscopic environment of a vertex will look like
\begin{figure}[H]
\includegraphics[scale=1]{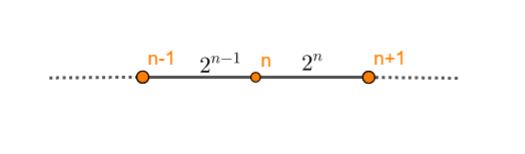}
\end{figure}
where the numbers on the edges are their respective conductances. Hence, $G$ is  transient. On the other hand, it is clear that this graph is not vertex-SIT, since any infinite simple path includes all the vertices in the $\N$ skeleton.

\end{document}